\documentclass[12pt]{article}

\usepackage{amssymb,amsthm,amsmath} %,amsrefs}
\usepackage[usenames,dvipsnames]{color}
\usepackage{soul}
\usepackage{longtable}
\allowdisplaybreaks
\usepackage{multirow}

\newcommand{\UMSF}{UMSF}

\newcommand{\comment}[1]{}
\definecolor{teal}{RGB}{0,128,128}
\definecolor{darkpurple}{RGB}{128,0,128}

\usepackage[titletoc,toc,title]{appendix}

\DeclareMathOperator{\lcm}{lcm}

\newtheorem{theorem}{Theorem}[section]

\newtheorem{lemma}[theorem]{Lemma}
\newtheorem{corollary}[theorem]{Corollary}

\theoremstyle{definition}
\newtheorem{definition}[theorem]{Definition}

\def \cF {{\cal F}}
\def \cG {{\cal G}}

\def \cT {{\cal T}}
\def \cU {{\cal U}}

\def \Z {\mathbb Z}

\title {On the minisymposium problem}
\author{
  P.\ Danziger 
  \thanks{Department of Mathematics, Toronto Metropolitan University, Toronto, ON M5B 2K3, Canada.
    E-mail: danziger@torontomu.ca}
  \and
  E.\ Mendelsohn
  \thanks{E-mail: mendelso@math.utoronto.ca}
  \and
  B.\ Stevens
  \thanks{School of Mathematics and Statistics - Carleton University, Ottawa, ON K1S 5B6, Canada. E-mail: brett@math.carleton.ca}
  \and 
  T.\ Traetta
  \thanks{DICATAM, Universit\`{a} degli Studi di Brescia, Via Branze 43, 25123 Brescia, Italy. E-mail: tommaso.traetta@unibs.it}
}
\date{}
\begin{document}
\maketitle

\begin{abstract}
The generalized Oberwolfach problem asks for a factorization of
  the complete graph $K_v$ into prescribed $2$-factors and at most 
  a $1$-factor. When all $2$-factors are pairwise isomorphic and $v$ is odd, we have the classic Oberwolfach problem, which
  was originally stated as a seating problem: given $v$ attendees at a conference with $t$ circular tables such that the $i$th table seats $a_i$ people and ${\sum_{i=1}^t a_i = v}$, find a seating arrangement over the 
$\frac{v-1}{2}$ days of the conference, so that every person sits next to each other person exactly once.
  
%  The Oberwolfach problem was originally stated as a seating problem: given $v$ attendees at a conference with $t$ circular tables such that the $i$th table seats $a_i$ people and ${\sum_{i=1}^t a_i = v}$, find a seating arrangement over the $r$ days of the conference, so that every person sits next to each other person exactly once.

  In this paper we introduce the related {\em minisymposium problem}, which requires
  a solution to the generalized Oberwolfach problem on $v$ vertices that contains a
  subsystem on $m$ vertices. That is, the decomposition restricted to the required $m$ vertices is a solution to the generalized Oberwolfach problem on $m$ vertices.  In the seating context above, the larger conference contains a minisymposium of $m$ participants, and we also require that pairs of these $m$ participants be seated next to each other for $\left\lfloor\frac{m-1}{2}\right\rfloor$ of the days.

  When the cycles are as long as possible, i.e.\ $v$, $m$ and $v-m$, a flexible method of Hilton and Johnson provides a solution. We use this result to provide further solutions when $v \equiv m \equiv 2 \pmod 4$ and all cycle lengths are even.
  In addition, we provide extensive results in the case where all cycle lengths are equal to $k$, solving all cases when $m\mid v$, except possibly when $k$ is odd and $v$ is even. 
  % In particular, we completely solve the case when all cycles are of length $m$.
\end{abstract}

\section{Introduction}

We assume that the reader is familiar with the fundamentals of graph theory and of design theory and refer them to \cite{West} and \cite{Handbook}, respectively.
In particular, a factor is a spanning subgraph and an $r$-factor is a factor which is $r$-regular, so in a 1-factor every vertex has degree one and a 2-factor is a disjoint union of cycles. 
Given a collection of factors, $\cal F$, an $\cal F$-factorization of a graph $G$ is a decomposition of the edges of $G$ into subgraphs, each of which is isomorphic to some $F\in {\cal F}$.
If ${\cal F} = \{F\}$ we speak of an $F$-factorization.

We use $K_n$ to denote the complete graph on $n$ vertices and $K_n^*$ to denote the graph $K_n$ when $n$ is odd and $K_n-I$, where $I$ is a 1-factor, when $n$ is even. 
Similarly, $K_{m,n}$ denotes the complete bipartite graph with parts of sizes $m$ and $n$. If the parts are $X$ and $Y$, respectively, we may also speak of $K_{X,Y}$. A 2-factor is called {\em uniform} if all of its constituent cycles are of the same length; it is called {\em Hamiltonian*} if its cycles have longest possible lengths given the requirements of the factorization. If a 2-factor, $F$, consists entirely of cycles of a particular length, $k$ say, we refer to an $F$-factor and an $F$-factorization as a $C_k$-factor and $C_k$-factorization, respectively.
Given a graph $G$, we denote by $G[n]$ the {\em lexicographic product} of $G$ with the empty graph on $n$ vertices. Specifically, the vertex set of $G[n]$ is $V(G)\times \Z_n$ (where $\Z_n$ denotes the cyclic group of order $n$) and 
$(x,i)(y,j)\in E(G[n])$ if and only if $xy \in E(G)$, $i,j\in \Z_n$.

The well known Oberwolfach problem asks for a 2-factorization of $K_n^*$ into 2-factors all of which are isomorphic to a given 2-factor $F$.
A summary of results up until 2006 can be found in \cite[Section VI.12]{Handbook}, in particular the case of uniform factors has been solved \cite{Alspach_Haggkvist_85, ASSW,Hoffman_Schellenberg_91}.
\begin{theorem}[\cite{Alspach_Haggkvist_85, ASSW,Hoffman_Schellenberg_91}]
  \label{uniform ober}
  Given integers $v, k\geq 3$, there is a $C_k$-factorization of $K_v^*$ if and only if $k\mid v$, except that there is no $C_3$-factorization of $K_6^*$ or $K_{12}^*$. 
\end{theorem}
% 1, 2, 16]
The case when all cycles in $F$ have even length has been completely solved in \cite{BryantDanziger}. The case with exactly two cycles is solved in \cite{Traetta_13}.
The case of the complete graph $K_\aleph$, where $\aleph$ is any infinite cardinal has been completely solved in \cite{Costa20}.
In the related Hamilton-Waterloo problem two 2-factors $F_1$ and $F_2$ are specified and we are asked for a factorization of $K_n^*$ into a given number of each of the factors. There has been much recent progress in this problem, see \cite{AsplundEtAl, BonviciniBuratti, BryantDanzigerDean, BDT1, BDT3, BDT2, DanzigerQuattrocchiStevens, KeranenOzkan, KeranenPastine,LeiShen,OdabasiOzkan,WangCao,WangChenCao}.
More generally, in the generalized Oberwolfach problem we are given a set of 2-factors $F_1, F_2, \ldots, F_t$ of $K_n$ and positive integers $\alpha_1, \alpha_2,\ldots,\alpha_t$, 
where $\sum \alpha_i = \lfloor \frac{n-1}{2} \rfloor$,
and are asked for a factorization of $K_n^*$ which contains $\alpha_i$ copies of the 2-factor $F_i$, see \cite{BryantDanziger, CaElKhoVan04, ElTipVan02}. A major recent development gives a non-constructive asymptotic existence result for the generalized Oberwolfach problem \cite{MR4012874}.

Other graphs have also been considered. In particular, Liu has shown the following for the complete multipartite graph.
\begin{theorem}[\cite{Liu00, Liu03}]
  \label{Liu} Let $k, t$ and $u$ be positive integers with $k\geq 3$.
  There exists a $C_k$-factorization of $K_t[u]$ if and only if $k\mid tu$, $(t-1)u$ is even, 
  further $k$ is even if $t = 2$, and
  $(k, t, u) \not\in \{(3, 3, 2), (3, 6, 2), (3, 3, 6), (6, 2, 6)\}$.
\end{theorem}

% The case where, very little.
Originally the Oberwolfach problem was stated as a seating problem:
\begin{quote}
  Given an odd number $v$ of attendees at a conference with $t$ circular tables such that the $i$th table seats $a_i$ people and ${\sum_{i=1}^t a_i = v}$, find a seating arrangement over the $\frac{v-1}{2}$ days of the conference, so that every person sits next to each other person exactly once.
\end{quote}

In this paper we introduce the related {\em minisymposium problem}. In this case we require a solution to the generalized Oberwolfach problem on $v$ vertices such that its restriction to a subset of $m$ vertices constitutes a solution to
the generalized Oberwolfach problem on $m$ vertices. Another way of considering the problem asks
% More simply stated, we are asking 
for a solution to the generalized Oberwolfach problem on $v$ vertices which contains a subsystem on $m$ vertices.
In the seating context above, the larger conference contains a minisymposium of $m$ participants, and we also require that pairs of these $m$ participants be seated next to each other for $\left\lfloor\frac{m-1}{2}\right\rfloor$ of the days.
A similar problem has been considered, for example, in \cite{BZ}
for whist tournaments.

Section~\ref{prelim} gives the formal definition of a minisymposium factorization and some
necessary conditions for its existence, as well as introduces some special cases. In Section~\ref{sec_ham_bi} we show how to use a flexible theorem by Hilton and Johnson \cite{MR1865547} to solve the
case of Hamiltonian* 2-factors (where the cycles are as long as possible). The same section considers the case where all cycles are of even length and $v \equiv m \equiv 2 \pmod 4$.
Section~\ref{Uniform Factors} considers the uniform case, where all cycles have the same length. We completely solve the case where all cycles are of length $m$ when $(v -1)m$ is even. In Section~\ref{sec_multiple_subsystems} we discuss and give some preliminary results on factorizations that contain more than one subsystem.  We provide some concluding remarks in the final section.

\section{Preliminaries}
% \section{Definitions and Necessary Conditions}
\label{prelim}

We begin by giving a formal definition of a minisymposium factorization. The
minisymposium problem is equivalent to the original Oberwolfach problem when
$v = m$. Hence we will generally assume that $v > m$.
% \pnote{Should we require $v>m$? This makes the statements of some of the Theorems later on easier (e.g. 2.4 5.?) The case $v=m$ is the original Oberwolfach problem.}
\begin{definition} \label{basic_def}
  Given positive integers $v$ and $m$ with $v\geq m$, let 
  $$\cF = \left\{F_i:\; i=1,\ldots, \left\lfloor\frac{v-1}{2}\right\rfloor-\left\lfloor\frac{m-1}{2}\right\rfloor\right\},$$ 
  be a collection of 2-factors on $v$ vertices and let 
  $$\cG = \left\{ (T_i, U_i) : i = 1,\ldots, \left\lfloor\frac{m-1}{2}\right\rfloor\right\},$$
  where the $T_i$ are 2-factors on $m$ vertices and the $U_i$ are 2-factors on $v-m$ vertices. 
  We define a {\em minisymposium factorization} MSF$(\cF, \cG)$ as a factorization of $K_v^*$ into 2-factors $F\in \cF$ and $G_i = T_i\cup U_i$, where $(T_i,U_i)\in \cG$, such that $\cT = \{T_i:\;  i = 1,\ldots, \left\lfloor\frac{m-1}{2}\right\rfloor\}$  is a factorization of a subgraph of $K_v^*$ isomorphic to $K_m^*$.
\end{definition}

Note that with the notation MSF$(\cF,\cG)$, we assume that the parameters $v$, $m$, $\cT$ and $\cU = \{U_i:\;   i = 1,\ldots, \left\lfloor\frac{m-1}{2}\right\rfloor\}$ are defined implicitly.
We may also use the notation MSF$(\cF,(\cT,\cU))$, if we wish to explicitly refer to the factorizations $\cT$ and $\cU$.

An MSF$(\cF,\cG)$ can be thought of as a 2-factorization of $K_v^*$ with a subsystem of size $m$.
When $m = 1$ or $2$, this is just a factorization of $K_v^*$ into 2-factors in $\cF$, which is equivalent to a solution of 
the 
generalized Oberwolfach problem, and so we will assume $m \geq 3$. Similarly, when $m = v$, this is a factorization into the $T_i$, so we assume that $m< v$.

Removing the subsystem, we can talk about a 2-factorization of $K_v^*$ with a ``hole'' of size $m$.
However, care must be taken when either $v$ or $m$ is even as the placement of the various 1-factors must 
be
considered, as noted below.

We note that the size of $\cF$ is
\[
  \left\lfloor\frac{v-1}{2}\right\rfloor-\left\lfloor\frac{m-1}{2}\right\rfloor= \left\{
    \begin{array}{cl}
      \frac{v-m}{2} & v \equiv m \pmod 2 \\
      \frac{v-m+1}{2} & v+1\equiv m \equiv 0 \pmod 2 \\
      \frac{v-m-1}{2} & v \equiv m+1 \equiv 0 \pmod 2. \\
    \end{array}
  \right.
\]

\begin{itemize}
\item
  In the case when both $v$ and $m$ are even, we are considering a factorization of $K_v^*=K_v-I$, where $I$ is a 1-factor, containing a factorization of a subgraph $G-J$ of $K_v^*$ where $G \cong K_m$ and $J$ is a 1-factor of $G$ 
contained in $I$. 

\item
  When $v$ is even and $m$ is odd, we are  considering a factorization of $K_v-I$, where $I$ is a 1-factor, containing a factorization of a subgraph $G \cong K_m$. Note that none of the edges of $I$ are contained in $G$. 

\item
  When $v$ is odd and $m$ is even, we are  considering a factorization of $K_v^*=K_v$ containing a factorization of $G-J$, where $G \cong K_m$ and $J$ is a 1-factor of $G$. We note that the edges of $J$ are not covered by factors in $\cT$, and hence must be covered by factors in $\cF$.

\item
  When both $v$ and $m$ are odd there is no 1-factor to consider.
\end{itemize}

Since in all cases $G \cong  K_m$, 
we will henceforth refer to it as {\em the} $K_m$. 
We note that
none of the edges of the $K_m$ are covered by $\cF$, except in the case of $m$ even and
$v$ odd; in this case, the edges of the 1-factor $J$ of the $K_m$ are covered by $\cF$. We use
this observation in the proofs of the lemmas below, where for a given 2-factor $F$,
we define $a_k(F)$ to be the number of cycles of length $k$ in $F$.
\begin{lemma}
  \label{General Necessary}
  For a given $v$ and $m$, if there is a minisymposium factorization, MSF$(\cF, \cG)$, then for each factor $F\in \cF$ using $c_F$ edges in the $K_m$,
  \begin{align}
    \label{v-m bound c_i}
    v-m & \geq -c_F +\sum_{i=3}^v a_i(F)\left\lceil\frac{i}{2}\right\rceil, \\
    \label{m bound c_i}
    m & \leq c_F + \sum_{i=3}^v a_i(F)\left\lfloor\frac{i}{2}\right\rfloor.
  \end{align}

  In particular, if $v$ is even, or $m$ is odd, all of the $c_F=0$ and therefore,
  \begin{equation}
    \label{ratio bound}
    \frac{\sum_{i=3}^v a_i(F)\left\lceil\frac{i}{2}\right\rceil}{\sum_{i=3}^v a_i(F)\left\lfloor\frac{i}{2}\right\rfloor} \leq \frac{v-m}{m}.
  \end{equation}

\end{lemma}
\begin{proof}
  We first deal with the case when $v$ is even, or $m$ is odd. In this case none of the edges of the $K_m$ appear in any $F\in \cF$.
  Therefore, for each $F \in \cF$, at most $\left\lfloor\frac{i}{2}\right\rfloor$ vertices of any cycle of length $i$ in $F$ are inside the $K_m$,
  hence 
  \begin{equation}
    \label{m bound}
    m \leq \sum_{i=3}^v a_i(F)\left\lfloor\frac{i}{2}\right\rfloor.
  \end{equation}
  Similarly, for any cycle of length $i$ in $F$, at least $\lceil\frac{i}{2}\rceil$ vertices of the cycle are not in the $K_{m}$, thus
  \begin{equation}
    \label{v-m bound}
    \sum_{i=3}^v a_i(F)\left\lceil\frac{i}{2}\right\rceil \leq v-m.
  \end{equation}
  Thus Inequality~\eqref{ratio bound} follows.

  Now, if $v$ is odd and $m$ is even the $\frac{m}{2}$ edges in the 1-factor $J$ of the $K_m$ must be used in factors from $\cF$.
  Suppose that $F\in \cF$ uses $c_F$ of these edges. Each edge of the $K_m$ used can increase the right hand side of Inequality~\eqref{m bound} by no more than one and decrease the left hand side of Inequality~\eqref{v-m bound} by no more than one.

  \comment{
    actually something stronger can be said if we think it is worth it. For a single cycle of length $i_j$ if it contains $c_{f,j}$ edges in the $K_m$ then the number of vertices it has in the $K_m$ is bounded above by $\lfloor \frac{i_j+c_{f,j}}{2} \rfloor$.
    If $i_j$ is odd then $\lfloor \frac{i_j+c_{f,j}}{2} \rfloor = \lfloor \frac{i_j}{2} \rfloor + \lceil \frac{c_{f,j}}{2} \rceil$. If $i_j$ is even then $\lfloor \frac{i_j+c_{f,j}}{2} \rfloor = \lfloor \frac{i_j}{2} \rfloor + \lfloor \frac{c_{f,j}}{2} \rfloor$.

    So if $c_f$ is larger than the number of odd cycles in $F$, $no_{F}$, then Inequality~(\ref{m bound c_i}) can be tightened to something like
    \[
      m \leq \sum_{i=3}^v \left(a_i(F)\left\lfloor\frac{i}{2}\right\rfloor\right)+no_{F} + \lfloor \frac{c_F-no_{F}}{2} \rfloor.
    \]

    However, when $F$ has only one odd cycle, for this bound to really constrain things would require that $v < 2m-1$ which never happens so We probably can ignore this.
    
  }

\end{proof}

\begin{theorem}
  \label{simple necessary condition}
  For a given $v$ and $m$, if there is a minisymposium factorization, MSF$(\cF, \cG)$, then $v\geq 2m$, unless $v$ is odd and $m$ is even, in which case $v\geq 2m-1$.
\end{theorem}
\begin{proof}
  When $v$ is even or $m$ is odd, the left hand side of Inequality~\eqref{ratio bound} is at least 1 and therefore $v\geq 2m$.

  When $v$ is odd and $m$ is even, $\sum_{F\in\cF} c_F = \frac{m}{2}$. Since $v$ is odd, each $F\in \cF$ must contain at least one odd cycle, therefore $\sum_{i=3}^v a_i(F)\lceil\frac{i}{2}\rceil \geq \lceil\frac{v}{2}\rceil = \frac{v+1}{2}$. Also note that the number of factors $F\in\cF$ is $\frac{v-m+1}{2}$.
  Summing Inequality~\eqref{v-m bound c_i} over all of the $F\in \cF$ twice gives
  \begin{align*}
    (v-m+1)(v-m) & \geq 2\sum_{F\in\cF}
                   \left(-c_F + \sum_{i=3}^v a_i(F)\left\lceil\frac{i}{2}\right\rceil\right)  \\
                 & = -m + 2\sum_{F\in\cF} \sum_{i=3}^v a_i(F)\left\lceil\frac{i}{2}\right\rceil \\
                 & \geq -m + 2\sum_{F\in\cF} \frac{v+1}{2}  \\
                 & = (v-m+1)\frac{v+1}{2}-m
  \end{align*}
  \comment{
    % Here is a more detailed derivation.
    Rearranging gives
    \[
      \begin{array}{rcl}
        % (v-m+1)(v-m) \geq (v-m+1)\left(\frac{v+1}{2}\right)-m \\
        (v-m+1)(v-m) - (v-m+1)\left(\frac{v+1}{2}\right)-m &\geq & 0 \\
        % (v-m+1)[(v-m) - \left(\frac{v+1}{2}\right)]-m &\geq& 0 \\
        % (v-m+1)[(\frac{v-1}{2}-m)] -m \geq 0 \\
        % \frac{1}{2}(v-m+1)[(v-1)-2m) -m \geq 0 \\
        (v-m+1)(v-2m-1) - 2m &\geq & 0 \\
        % v^2 -2mv - v - vm +2m^2 + m + v - 2m - 1 - 2m \geq 0\\
        v^2 - 3mv - 3m + 2m^2 &\geq & 0\\
        % (v-(2m-1))(v-(m+1))\\
        % = v^2 - [(m+1)+(2m-1)]v +(2m-1)(m+1) \\
        % = v^2 - 3mv + 2m^2 + 3m  - 1  \\
        (v-(2m-1))(v-(m+1)) & \geq & 0
      \end{array}
    \]
  }

  Hence,
  \begin{align*}
    0 & \leq (v-m+1)(v-m)- (v-m+1)\frac{v+1}{2} + m  \\
      & =  \frac{(v-m+1)(v-2m-1) + 2m}{2} \\
      & =\frac{ (v - (2 m - 1)) (v - (m + 1))}{2}.
  \end{align*}
  When $v=m+1$, the factors in $\cU$ are required to be 2-regular graphs on a single vertex, which is not possible, so $v-(m+1) > 0$. Thus $v-(2m-1) \geq 0$ and the result follows.
\end{proof}
% correction for above
% 0 & \leq (v-m+1)(v-m)- (v-m+1)\left(\frac{v+1}{2}\right) + m  \\
% & \leq  \frac{(v-m+1)(v-2m-1) + 2m}{2} \\
% & =\frac{ (v - (2 m - 1)) (v - (m + 1))}{2}.

In the case where $v$ is odd and $m$ is even and $v = 2m-1$, we have that $\cU$ is a factorization of $K_{v-m}$.  So we can interchange both the roles of $m$ and $v-m$, as well as those of $\cT$ and $\cU$.
Thus, without loss of generality, we may assume that $v\geq 2m$ in all cases.

There are two cases of initial special interest. Firstly, the case of {\em uniform} cycle lengths (when all cycles in a factor are of the same length), which we consider in detail in Section~\ref{Uniform Factors}. Secondly, the case where the cycles are as long as possible, which in correspondence with the definition of $K_n^*$ and the Hamiltonian-like nature of such factorizations we will call {\em Hamiltonian*} factorizations. We formally define {\em Hamiltonian*} factorizations in Section~\ref{sec_ham_bi}.  There we show that a method of Hilton and Johnson completely settles their existence.

An MSF$(\cF,\cG)$ in which all of the factors in $\cF$, $\cT$ and $\cU$  are uniform with the same cycle length $k$ is called {\em uniform} and 
we
refer to it as a \UMSF$(v,m,k)$. 
In this case we have the following necessary conditions.
\begin{theorem}
  \label{k necessary condition}
  If $v>m$ and a \UMSF$(v,m,k)$ exists, then $k\geq 3$, $k\mid m$ and $k\mid v$. Furthermore, 
  \begin{itemize}
  \item
    if $k$ is even, then $v\geq 2m$;
  \item 
    if $k$ is odd, then $v \geq \frac{2mk}{k-1}$.
  \end{itemize}
\end{theorem}
\begin{proof}
  Since we are forming 2-factors with cycles of length $k$, we require $k\geq 3$.
  The divisibility conditions follow directly from the requirement for a factorization of $K_v^*$ and $K_m^*$ into $k$-cycles.
  If $k$ is even, then $v$ is even, since it is a multiple of $k$, and  Theorem~\ref{simple necessary condition} gives $v\geq 2m$.

  If $k$ is odd, we note that for a $C_k$-factor $F\in\cF$, we have $a_i(F) = \frac{v}{k}$ when $i=k$ and 0 otherwise, $\lfloor\frac{k}{2}\rfloor = \frac{k-1}{2}$ and $\lceil\frac{k}{2}\rceil = \frac{k+1}{2}$. Thus, when $v$ is even or $m$ is odd, Inequality~\eqref{ratio bound} implies that $v \geq \frac{2mk}{k-1}$ and the result follows.

  This leaves the case when $v$ and $k$ are odd and $m$ is even. We sum Inequality~\eqref{v-m bound c_i} over all the $F\in\cF$ to obtain
  \[
    \begin{array}{rcl}
      % v-m & \geq & \sum_{i=3}^v \left(a_i(F)\left\lceil\frac{i}{2}\right\rceil\right)-c_F \\
      \sum_{F\in\cF} (v-m) & \geq & \sum_{F\in\cF}\left(-c_F +  \sum_{i=3}^v a_i(F)\left\lceil\frac{i}{2}\right\rceil\right) \\
      % \frac{1}{2}(v-m+1)(v-m) & \geq & \left(\sum_{F\in\cF} \sum_{i=3}^v  \left(a_i(F) \left \lceil \frac{i{2} \right\rceil \right) \right) - m/2 \\
      \frac{1}{2}(v-m+1)(v-m) & \geq & -m/2 + \sum_{F\in\cF} \frac{v}{k}\frac{k+1}{2} \\
      % (v-m+1)(v-m) & \geq &  (v-m+1)\frac{v}{k}(k+1)/2-m \\
      2k(v-m+1)(v-m) & \geq &  v(v-m+1)(k+1)-2mk. \\
    \end{array}
  \]
  \comment{
    % Here is the derivation of the quadratic.
    So
    \[
      \begin{array}{rcl}
        2k(v-m+1)(v-m) - v(v-m+1)(k+1)+2mk & \geq &  0 \\
        2k(v^2 -2mv + m^2 + v - m) - (k+1)(v^2 - mv + v) + 2mk & \geq &  0 \\
        (k-1)v^2 + 2k(-2m + 1)v + (k+1)(m - 1)v + 2km^2 - 2km + 2mk & \geq & 0 \\
        (k-1)v^2 + (k(-4m + 2 + m - 1) + m - 1)v + 2km^2  & \geq & 0 \\
        (k-1)v^2 + (k(-3m  + 1) +  m - 1)v + 2km^2 & \geq & 0 \\
        (k-1)v^2 + (-3km  + k +  m - 1)v + 2km^2  & \geq & 0, \\
      \end{array}
    \]
  }
  Rearranging and expanding in $v$ gives
  \begin{equation}
    \label{v quadratic}
    (k-1)v^2 + (-3km+m+k-1)v + 2km^2 \geq 0.
  \end{equation}
  Let $$f(v)= (k-1)v^2 + (-3km+m+k-1)v + 2km^2.$$
  By Theorem~\ref{simple necessary condition}, $v \geq 2m-1$, but
  $$f\left(2m-1\right) = m(1+k - 2 m) < 0,$$
  so $v$ is at least as large as the larger of the two roots of $f$. Now
  $$
  f\left(\frac{2mk}{k-1} -2\right) = -2(m+1-k) < 0, \mbox{ and } f\left(\frac{2mk}{k-1} -1\right) = (k-1)m > 0.
  $$
  Thus $f$ has its larger root between $\frac{2mk}{k-1} -2$ and $\frac{2mk}{k-1} -1$.

  It is left to check that $v\neq \lfloor\frac{2mk}{k-1}-1\rfloor, \lceil\frac{2mk}{k-1}-1\rceil$.
  Recalling that $k\mid v$, if $v=\lfloor\frac{2mk}{k-1}-1\rfloor$ or $v=\lceil\frac{2mk}{k-1}-1\rceil$, 
  then there exist a rational number 
  $0\leq \epsilon < 1$ and a positive integer $a$ such that
  $$v = \frac{2mk}{k-1}-1 \pm \epsilon = ak.$$
  Multiplying both sides by $\frac{k-1}{2k}$ and rearranging we have that
  $$m- a\frac{k-1}{2} =  (1 \pm \epsilon)\frac{k-1}{2k}.$$
  Since $k$ is odd, the left side is an integer. 
  However, $0< (1 \pm \epsilon)\frac{k-1}{2k} <1$, a contradiction.  We conclude that $v\geq \frac{2mk}{k-1}$.
\end{proof}

One interesting case in light of these necessary conditions is when $m=k$, i.e.\ a \UMSF$(v,k,k)$. For these parameters, since $k < v$ and $k \mid v$, we must have $v \geq 2k$, so the necessary conditions in Theorem~\ref{simple necessary condition} are satisfied. We consider these types of factorizations in Section~\ref{Uniform Factors}.

\section{Hamiltonian* and bipartite factors}\label{sec_ham_bi}

Considering non-uniform factors, an obvious case to consider is a {\em Hamiltonian* minisymposium factorization}, which is one in which the cycles have the longest possible lengths. Specifically, the factors in $\cF$ are all $v$-cycles, factors in $\cT$ are all $m$-cycles and the factors in $\cU$ are all $(v-m)$-cycles. Such an MSF$(\cF,\cG)$ is denoted by HMSF$(v,m)$. We sometimes refer to the cycles in $\cT$ and $\cU$ as `short' cycles. 
Because of the lengths of these cycles there are no further necessary conditions beyond those of Theorem~\ref{simple necessary condition}.

In a paper on the Oberwolfach problem, Hilton and Johnson prove the following theorem on a flexible construction technique.
\begin{theorem}[\cite{MR1865547}] \label{hilton_thm}
  Let $m$ and $n$ be integers, $1 \leq m < n$. Let $(s_1,\ldots, s_t )$, $s_i \in {1, 2}$,
  $1 \leq i \leq t$, be a composition of $n-1$. Let $K_m$ be edge coloured with $t$ colours $c_1,\ldots, c_t$.
  Let $f_i$ be the number of edges coloured $c_i$ and $K_m(c_i)$ be the $i$th colour class. This colouring can be extended to an edge-colouring of $K_n$ in which the colour class $K_n(c_i)$ is an $s_i$-factor, $1 \leq i \leq t$, and when $s_i = 2$, $K_n(c_i)$ contains exactly one more cycle than $K_m(c_i)$ if and only if for all $i = 1,2, \ldots, t$:
  \begin{align*}
    f_i \geq s_i( m-n/2),  \\
    s_in \mbox{ is even},\\
    \Delta(K_m(c_i)) \leq s_i.\\
  \end{align*}
\end{theorem}

This theorem is sufficient to provide a solution to the Hamiltonian* minisymposium factorization.
\begin{corollary}
  \label{HMS exist}
  An HMSF$(v,m)$ exists if and only if $m\geq 3$, $v \geq 2m-1$ in case $m$ is even, and $v \geq 2m$ in case $m$ is odd.
\end{corollary}
\begin{proof}
  The given conditions are necessary by Theorem~\ref{simple necessary condition}.
  To prove sufficiency, we will define an edge colouring of the $K_m$ from a decomposition of the $K_m$ into Hamiltonian cycles and possibly a single 1-factor using Theorem~\ref{uniform ober}.  If $m$ is odd, this defines edge colours $c_i$, $1 \leq i \leq (m-1)/2$.  If $v$ is also odd, extend this to a $(v-1)/2$-edge colouring of the $K_m$ by including $(v-m)/2$ empty colour classes $c_i$, $(m+1)/2 \leq i \leq (v-1)/2$. Let $s_i = 2$ for all $1 \leq i \leq (v-1)/2$. If $v$ is even, extend the colouring to a $v/2$-edge colouring of the $K_m$ by adding empty colour classes.  Let $s_i = 2$ for all $1 \leq i < v/2$ and $s_{v/2} = 1$.  In both cases, it can be verified that Theorem~\ref{hilton_thm} now gives an HMSF$(v,m)$ as desired.

  If $m$ is even, define $m/2$ edge colour classes of the $K_m$ from a decomposition into Hamiltonian cycles and one 1-factor.  Let $c_1$ be the colour class of the 1-factor.  If $v$ is also even, extend this to a $v/2$-edge colouring by adding empty colour classes.  Let $s_1 = 1$ and $s_i = 2$ for $2 \leq i \leq v/2$.  If $v$ is odd, extend this to a $(v-1)/2$-edge colouring by adding empty colour classes.  Let $s_i = 2$ for $1 \leq i \leq v/2$.  In both cases, it can be verified that Theorem~\ref{hilton_thm} now gives an HMSF$(v,m)$ as desired.
\end{proof}

Theorem~\ref{hilton_thm} is more than just an existence result; a recursive procedure can be extracted from the proof to algorithmically build the edge decompositions.  We have a more direct construction of all HMSF$(v,m)$ which uses difference methods and decompositions of Cayley graphs \cite{hamiltonian_paper}.

Theorem~\ref{hilton_thm} can be used much more generally to build minisymposium factorizations.  Essentially it shows that it is possible to extend any 2-factorization of the $K_m$ to one of $K_v$, provided that the necessary conditions hold, where the additional 2-factors are Hamiltonian, with an additional 1-factor when $v$ is even.

When all of the cycles of the factors in $\cF$, $\cU$ and $\cT$ are bipartite (i.e.\ contain only even cycles), 
we apply the Theorem of H\"{a}ggkvist \cite{Hag} (given below) to HMSF$(v,m)$ to give us a solution to the minisymposium problem when $v \equiv m\equiv 2 \pmod 4$.

\begin{theorem}[\cite{Hag}]
  \label{Hagg}
  If $F$ is a bipartite 2-regular graph of order $2n$, then there is a factorization of $C_n[2]$ into 2 isomorphic copies of $F$.
\end{theorem}

We note that in the case where the factors are bipartite, so all cycle lengths are even, $v$ and $m$ are both even and Theorem~\ref{simple necessary condition} gives $v\geq 2m$. 

\begin{theorem}
  \label{Hagg Result}
  If  $v\equiv m\equiv 2\pmod 4$, 
  and $\cF = \{ F_i:\; 1\leq i\leq (v-m)/2\}$, $\cU = \{U_j :\; 1\leq j\leq (m-2)/2 \}$ and $\cT=\{T_j:\; 1\leq j\leq (m-2)/2\}$ are sets of bipartite factors with $F_i=F_{i+1}$, $U_j=U_{j+1}$ and $T_j=T_{j+1}$ for every odd $i$ and $j$, then an MSF$(\cF, (\cT, \cU))$ exists
  if and only $v\geq 2m$.
\end{theorem}
\begin{proof}
  We note that if $v \equiv m\equiv 2\pmod 4$, the number of the $F_i$ is $(v-m)/2$ and the number of the $U_j$ and the $T_j$ is $(m-2)/2$, so the number of 
  the
  $F_i$, $U_j$ and $T_j$ are all even.
  We take an HMSF$(v/2,m/2)$, which exists by Corollary~\ref{HMS exist}, with factors $F_i'$ of order $v/2$, $U_i'$ of order $(v-m)/2$ and $T'_i$ of order $m/2$.
  We blow up each vertex by 2 and apply Theorem~\ref{Hagg} to factor each $F_i'[2]$ into 2 copies of $F_i$, each $U_j'[2]$ into 2 copies of $U_j$ and each $T_j'[2]$ into 2 copies of $T_j$.
\end{proof}

If $v \equiv 0 \pmod 4$ or $m \equiv 0 \pmod 4$, then $v/2$ or $m/2$ would be even and the HMSF$(v/2,m/2)$ would contain 1-factors either in $K_{v/2}$ or the $K_{m/2}$.  When a 1-factor is blown up as done in Theorem~\ref{Hagg}, it results in a
%set of 4-cycles 
$C_4$-factor, which prevents constructing the desired MSF unless $\cF$, $\cT$, and $\cU$ already contain this kind of factor.  

An immediate consequence of Theorem~\ref{Hagg Result} is the following relating to uniform factors.
\begin{corollary}
  \label{even cycles}
  If $k>3$, $k\equiv 2 \pmod 4$, $v \equiv m\equiv k \pmod {2k}$, then a \UMSF$(v,m,k)$ exists if and only if $v \geq 2m$.
\end{corollary}

\section{Uniform Factors}
\label{Uniform Factors}

In this section we consider the case of uniform factors, i.e.\ when all cycles are of the same length, $k$. 
We recall from Theorem~\ref{k necessary condition} that in order for a \UMSF$(v,m,k)$ to exist, we require that $k\geq 3$, which we will assume throughout this section. We also require $k\mid m$ and $k\mid v$.  Additionally, if $k$ is even, then $v\geq 2m$ and if $k$ is odd, then $v \geq \frac{2mk}{k-1}$.

Corollary~\ref{even cycles} gives a powerful result in the case when $k\equiv 2 \pmod 4$ and $v \equiv m\equiv k\pmod {2k}$.
The case where $k=3$ has been considered in \cite{Rees_Stinson_1, Rees_Stinson_2, Stinson} when $v$ and $m$ are both odd, and \cite{DRS_08, DRS_03_1, DRS_03_2, GR, TRS} when they are both even. However, the case when $m$ and $v$ have opposite parities appears to be completely open. We summarize these results in the following theorem.
\begin{theorem}[\cite{Rees_Stinson_2, DRS_08}]
  \label{k=3}
  If $v\equiv m \pmod 2$, there exists a \UMSF$(v, m, 3)$ if and only if 
  $v \geq 3m$,  
  $v \equiv m \equiv 0 \pmod 3$,
  and if $v,m$ are even, then $v,m > 12$.
\end{theorem}

We will find the following results useful.
A corollary of a result in \cite{BS81} yields the following.
\begin{theorem}[\cite{BS81}]
  \label{BS81}
  If $G$ is a Hamiltonian decomposable graph, then $G[n]$ is also Hamiltonian decomposable. In particular, $C_m[n]$ has a $C_{mn}$-factorization for every $m\geq 3$.
\end{theorem}

Piotrowski \cite{Piotrowski} has shown the following result for $m \geq 4$. The case $m=3$ is covered by Theorem~\ref{Liu}.
\begin{theorem}[\cite{Piotrowski}]
  \label{C_m || C_m[n]}
  There exists a $C_m$-factorization of $C_m[n]$, except if $n=2$ and $m$ is odd, or when $(m,n) =(3,6)$.
\end{theorem}
Piotrowski \cite{Piotrowski} has also shown the following result.
\begin{theorem}[\cite{Piotrowski}]
  \label{bipartite}
  Let F be a bipartite 2-regular graph of order $2n$. The complete bipartite graph $K_2[n]$ has an $F$-factorization if and only if $n$ is even, except when $n = 6$ and $F$ consists of two 6-cycles.
\end{theorem}

We now give some recursive constructions for uniform minisymposium factorizations. 
\begin{theorem}
  \label{mt construction}
  Let $m\geq k\geq 3$ and $t\geq 2$ be integers.  
  If $(t-1)m$ is even and $k\mid m$, then there is a \UMSF$(mt, m, k)$, except that there is no \UMSF$(6t, 6, 3)$ \UMSF$(12t,12,3)$, \UMSF$(12, 6, 6)$, or \UMSF$(2m, m, k)$ when $k$ is odd.
\end{theorem}
\begin{proof}
  The non-existence of a \UMSF$(6t, 6, 3)$ and \UMSF$(12t, 12, 3)$ are covered by Theorem~\ref{k=3}.
  Since a \UMSF$(2m, m, k)$ is equivalent to a $C_k$-factor\-ization of the complete bipartite graph $K_{2}[m]$, 
  it clearly does not exist when the cycle length $k$ is odd, or when $k=m=6$ by Theorem~\ref{bipartite}. 
  In all remaining cases, the following conditions simultaneously hold:
  \begin{enumerate}
  \item $(m,k)\not\in\{(6, 3), (12, 3)\}$, 
  \item $(t,m,k)\neq(2,6,6)$, 
  \item if $k$ is odd, then $t>2$.
  \end{enumerate}
  The assumptions of Theorems~\ref{uniform ober} and \ref{Liu} are
  then
   satisfied. Hence there is a $C_k$-factorization of $K_t[m]$ and a $C_k$-factorization of $K^*_m$, which we use to fill in the parts of size $m$ in $K_t[m]$. This completes the proof.
\end{proof}

Considering the necessary conditions in Theorem~\ref{k necessary condition}, we get the following corollaries.
\begin{corollary}
  \label{m|v}
  Suppose that either $k$ is even or $v$ is odd, and $m\mid v$. Then there exists a \UMSF$(v,m,k)$ if and only if $k \mid m$, $v\geq 2m$ when $k$ is even, and $v \geq 3m$ when $k$ is odd, except that \UMSF$(v,6,3)$, \UMSF$(v, 12, 3)$ and \UMSF$(12,6,6)$ do not exist.
\end{corollary}
\begin{proof}
  Taking $v=mt$, Theorem~\ref{k necessary condition} gives the necessary conditions $k \mid m$ and  $v\geq 2m$ when $k$ is even. When $k$ is odd, the necessary condition from Theorem~\ref{k necessary condition} is $v \geq \frac{2mk}{k-1}$, but since $m\mid v$, this implies $k \geq 3m$.
  Given the conditions of $k$ and $v$, the sufficiency comes from Theorem~\ref{mt construction}.
\end{proof}

We note that if $m\mid v$ this corollary completely solves all cases except when $k$ is odd and $v$ is even.
One case of particular interest is when $k=m$, in this case $m\mid v$ is necessary.
\begin{corollary}
  \label{k=m}
  Let $m(t-1)$ be even.  Then a \UMSF$(tm,m,m)$ exists if and only if $t \geq 2$ when $m$ is even, $t\geq 3$ when $m$ is odd, and $(t,m) \neq (2,6)$.
\end{corollary}

The previous results all require $m~\mid~v$, however the next theorem allows us to recursively construct solutions to cases where $m$ does not divide $v$. 

\begin{theorem}
  \label{recursive}
  Assume there is a \UMSF$(v,m,k)$ and let $t\geq 1$. Then there exists a \UMSF$(vtk, mtk, \ell)$, with $\ell\in\{k, kt\}$, in each of the following cases:
  \begin{enumerate}
  \item $v$ and $m$ have the same parity;
  \item $v$ and $t$ are even, $\ell=tk$, and $m$ and $k$ are both odd, except possibly when 
    $(k, t) = (3,2)$.
  \end{enumerate}
\end{theorem}
\begin{proof} 
  Letting $w\in\{m,v\}$, we factorize $K^*_{wtk}$ into $\Gamma_w = K^*_w[tk]$ and $\overline{\Gamma}_w = K^*_{wtk}- \Gamma_w$. 
  Note that $\overline{\Gamma}_w$ is the vertex disjoint union of
  \begin{enumerate}
  \item $w$ copies of $K^*_{tk}$ when $w$ is odd, or
  \item $w/2$ copies of $K^*_{2tk}$ when $w$ is even.
  \end{enumerate}
  Without loss of generality, we can assume that $\Gamma_m \subseteq \Gamma_v$ and $\overline{\Gamma}_m \subseteq \overline{\Gamma}_v$, except when $v$ is odd and $m$ is even. 
  In this case the components of $\overline{\Gamma}_m$ are copies of $K^*_{2tk}$, while those of $\overline{\Gamma}_v$ are isomorphic to $K^*_{tk}$, therefore $\overline{\Gamma}_m \subseteq \overline{\Gamma}_v$ cannot hold. 
  We proceed by constructing 
  \begin{enumerate}
    \renewcommand{\theenumi}{$\alph{enumi}$}
    \renewcommand{\labelenumi}{(\theenumi)}
  \item\label{a} a $C_\ell$-factorization of $\Gamma_v$
    containing a $C_\ell$-factorization of $\Gamma_m$, and
  \item\label{b} a $C_\ell$-factorization of $\overline{\Gamma}_v$ 
    containing a $C_\ell$-factorization of $\overline{\Gamma}_m$,
  \end{enumerate}
  which together will provide  the desired
  $\UMSF(vtk, mtk, \ell)$.

  We blow up each vertex of the \UMSF$(v,m,k)$ by $tk$, to obtain a $C_k[tk]$-factorization of $\Gamma_v$ containing a $C_k[tk]$-factorization of $\Gamma_m$. To construct (\ref{a}) it is therefore enough to factorize $C_k[tk]$ into $C_\ell$-factors, $\ell\in\{k, tk\}$.
  By Theorem~\ref{C_m || C_m[n]} there is a $C_k$-factorization of $C_k[tk]$, except when $(t,k)=(2,3)$. In this case, the desired  \UMSF$(6v,6m,3)$ exists by Theorem \ref{k=3}.
  Considering that $C_k[tk] = C_k[t][k]$, by Theorem~\ref{BS81} there exists a $C_{kt}$-factorization of $C_k[t]$ which we blow up by $k$ to obtain  $C_{kt}[k]$-factorization of $C_k[tk]$. 
  By Theorem~\ref{C_m || C_m[n]}, each $C_{kt}[k]$-factor can be further decomposed into $C_{kt}$-factors yielding a  $C_{kt}$-factorization of $C_k[tk]$. 

  It is left to construct (\ref{b}). 
  If $m$ and $v$ have the same parity, the components of $\overline{\Gamma}_m$ and $\overline{\Gamma}_v$ are pairwise isomorphic: they are copies of $K^*_{tk}$ or $K^*_{2tk}$. It is then enough to build a $C_{\ell}$-factorization of $K^*_{tk}$ and $K^*_{2tk}$ for $\ell\in\{k, kt\}$. They exist by Theorem~\ref{uniform ober} except when $\ell = k = 3$ and one of the following two conditions hold,
  \begin{enumerate}
  \item $mv$ is odd and $t\in \{2, 4\}$, or
  \item $m$ and $v$ are even, and $t\in \{1, 2\}$. 
  \end{enumerate}
  In each of these cases, the existence of the desired $\UMSF(3vt, 3mt, 3)$ is guaranteed by Theorem \ref{k=3}.

  If $v$ and $t$ are even, $\ell=tk$, and both $m$ and $k$ are odd, the components of $\overline{\Gamma}_m$ are isomorphic to $K^*_{tk}$, while those of $\overline{\Gamma}_v$ are isomorphic to $K^*_{2tk}$. 
  Since we can factorize $K^*_{2tk}$ into $K_2[tk]$ and two copies of $K^*_{tk}$, it is enough to decompose both $K_2[tk]$ and $K^*_{tk}$ into $C_{tk}$-factors.
  These factorizations exist by Theorem~\ref{bipartite}  and Theorem~\ref{uniform ober}, respectively, except possibly when $(k, t) = (3,2)$.
\end{proof}

We may now use the result on triples (Theorem~\ref{k=3}) to obtain the following.
\begin{corollary}
  \label{recurse triples}
  Let $v\equiv m \equiv 0,3\pmod 6$, with  $v \geq 3m$ and $m\not\in\{0,6,12\}$. Then 
  there exists a \UMSF$(3tv, 3tm,3t)$ for all $t>0$.
\end{corollary}

Additionally, we may use Theorem~\ref{even cycles} to obtain the following result.
\begin{corollary}
  \label{recurse bipartite}
  Let $3<k$, $k \equiv 2 \pmod 4$, $v \equiv m\equiv k\pmod {2k}$ and $v \geq 2m$.  Then there exists a \UMSF$(vtk, mtk,k)$ and a \UMSF$(vtk, mtk,tk)$ for all $t>0$.
\end{corollary}

We note that the above result can be used to obtain \UMSF's with cycle length, subsystem size or number of vertices congruent to $0 \pmod 4$ by taking $t$ even.  However, in all cases, the number of vertices and subsystem size will be divisible by $k^2$.

\section{Multiple Subsystems}\label{sec_multiple_subsystems}

A natural question to ask is if a system can have multiple subsystems. 
In general, it seems likely to be hard 
to navigate through the lattice of subsystems and all the possible ways the subsystems can be distributed across the main system.  However, when the subsystems are disjoint, have small common intersections or are nested, the problem is more tractable.  We give some preliminary results in the next three subsections.

\subsection{Disjoint Subsystems}

In the uniform case, the flexibility of Theorem~\ref{Liu} allows us to create a large number of disjoint subsystem. We refer to a factorization of $K_{v}^*$ into $k$-cycles with subsystems on disjoint vertex sets of sizes $m_j$ for $1 \leq j \leq n$ 
as a \UMSF$(v,\{m_j\},k)$.

\begin{lemma}
  \label{disjoint_holes_minimal}
  Let $k \mid m_j$ for $1 \leq j \leq n$. Let $m$ be an integer such that there is a \UMSF$(m,m_j,k)$ for each $1 \leq j \leq n$.  Then there exists a \UMSF$(ms,\{m_j\},k)$ for 
  all
  $s \geq \max\{2,n\}$ if $k$ is even, 
  and for all $s \geq \max\{3,n\}$ such that $(s-1)m$ is even if $k$ is odd, except when $(k, s, m) =(6, 2, 6)$.
\end{lemma}
\begin{proof}
  Theorem~\ref{Liu} guarantees the existence of a $C_k$-factorization of $K_s[m]$.  
  For each $1 \leq j \leq n$, place a copy of the ingredient \UMSF$(m,m_j,k)$ on the $j$th part of size $m$ of  $K_s[m]$, and any $C_k$-factorization of $K_m^*$ on each of the remaining parts.
  The definite exception $(k, s, m) =(6, 2, 6)$ follows from the non-existence of
  a \UMSF$(12, 6, 6)$
  (see Theorem \ref{mt construction}).
\end{proof}

As with the uniform factorizations containing a single subsystem in this paper, the easiest case is when $m_j \mid m$ for all $1 \leq j \leq n$ and either $k$ is even or $m$ is odd.

\begin{corollary}
  \label{disjoint_holes_thm}
  Let $m = \lcm\{m_j:\; j = 1,\ldots,n\}$, and assume the following conditions are all satisfied:
  \begin{enumerate}
  \item $k \mid m_j$ for all $j \in \{1,\ldots,n\}$ and $m \mid v$;
  \item $k(m-1)$ is even;
  \item if $k=3$, then $m_j \not\in \{6,12\}$ for all $j$;
  \item if $(k,m)=(6,12)$, then $m_j \neq 6$ for all $j$;
  \item $(v,m,k) \neq (12,6,6)$;
  \item $v/m \geq \max\{2,n\}$ if $k$ is even;
  \item $v/m \geq \max\{3,n\}$ and $v$ is odd if $k$ is odd.
  \end{enumerate}
  Then there exists a \UMSF$(v,\{m_j\},k)$.
\end{corollary}
\begin{proof}
  Since each $m_j$ is a divisor of $m$ and conditions $1-4$ hold, we can apply either Theorem \ref{uniform ober} or Corollary \ref{m|v}, as needed, to ensure the existence of a \UMSF$(m,m_j,k)$ for every $j \in \{1,\ldots,n\}$.
  These are the ingredient designs needed in Lemma~\ref{disjoint_holes_minimal}, which can be applied in view of conditions $5-7$.
\end{proof}

We note that for any fixed multiset of $m_j$, this corollary constructs \UMSF$(v,\{m_j\},k)$ for all but a finite number of $v$ permitted by the necessary conditions when $m_j \mid v$ and either $k$ is even or $m$ is odd. 

\subsection{Scattered Subsystems}

The proof of Lemma \ref{disjoint_holes_minimal} builds systems whose factors 
intersect either all of the subsystems or none of them. A balancing of the sizes of these intersections could be an interesting property. For instance,
we could ask for systems whose factors do not intersect more than one subsystem.
In other words,
we ask for a $C_k$-factorization $\mathcal{F}$ of $K_{v}^*$ that contains 
$n$ subsystems 
% $\mathcal{S}_1, \mathcal{S}_2, \ldots, \mathcal{S}_n$ 
of sizes $m_1, m_2, \ldots, m_n$, such that no two factors of any of the subsystems
are contained in the same factor of $\mathcal{F}$. We denote such a factorization by \UMSF$(v, [m_j], k)$ and say that the subsystems are {\em scattered}.

Partial results in this direction could be easily obtained by making use
of cycle frames. We recall that a {\em $k$-cycle frame} ($k$-CF) of $K_{s}[m]$ is a decomposition of
$K_{s}[m]$ into holey $C_k$-factors; a {\em holey $C_k$-factor} is a vertex-disjoint union of $k$-cycles 
covering all vertices $K_{s}[m]$ except those belonging to one part.
The following result, proven in \cite{BCDT17}, provides necessary and sufficient conditions for the existence of $k$-cycle frames.
\begin{theorem}[\cite{BCDT17}]\label{CF} Let $m\geq 2$ and $k,s\geq 3$.
  There exists a $k$-cycle frame of $K_{s}[m]$ if and only if
  $m$ is even, $m(s-1)\equiv 0\pmod{k}$, $k$ is even when $s=3$, and
  $(k, m, s) \neq (6,6,3)$. 
\end{theorem}
By making use of Theorem \ref{CF}, we obtain the following.

\begin{lemma}\label{scattered_subsystems_minimal}
  Let $u\in\{1,2\}$.
  If there exists a \UMSF$(2m+u, m_j, k)$ for each $1 \leq j \leq n$, 
  then there exists a \UMSF$(2ms+u, [m_j],k)$ whenever $2s\equiv 2 \pmod{k}$ and
  $s \geq n$.
\end{lemma}
\begin{proof} Let $n\geq 1$,
  $2s\equiv 2 \pmod{k}$ and $s \geq n$. 
  It follows that 
  $2m(s-1)\equiv 0\pmod{k}$, $k=4$ when $s=3$, and $(k, 2m, s) \neq (6,6,3)$. 
  Therefore, Theorem \ref{CF} guarantees the existence of a 
  $k$-cycle frame $\mathcal{F}$ of $K_{s}[2m]$. Let $P_i$ denote the $i$-th part
  of $K_{s}[2m]$, for $1\leq i\leq s$. Also, let 
  \[\mathcal{F}=\{F_{ij}:\; 1\leq i\leq s, 1\leq j\leq m\},
  \]
  where the $F_{ij}$s are the holey $C_{k}$-factors of $\mathcal{F}$ missing $P_i$, for
  $1\leq i\leq s$.
  By assumption, there is a \UMSF$(2m+u, m_j, k)$ on $P_i \cup \{\infty_1, \infty_u\}$,
  say $\mathcal{H}_i = \{H_{ij}:\;  1\leq j\leq m\}$. It follows that
  $\mathcal{F}^*=\{F_{ij} \cup H_{ij}:\; 1\leq i\leq s, 1\leq j\leq m\}$ is
  a $C_k$-factorization of $K_{2ms+u}$ with scattered subsystems of sizes $m_1, m_2, \ldots, m_n$.
  Indeed, the factors of the subsystems belong to the $H_{ij}$s, each of which belongs to exactly one factor of $\mathcal{F}^*$.
\end{proof}
In the \UMSF$(2ms+u, [m_j],k)$ constructed in the proof of Lemma~\ref{scattered_subsystems_minimal}, two subsystems may intersect in 0, 1, or 2 vertices, which are necessarily in the set $\{\infty_1,\infty_u\}$.
% Given subsystems of sizes $m_1, m_2, \ldots, m_n$, they may pairwise intersect in 0, 1 or 2 vertices. By placing the subsystem of size $m_j$ on 0, 1 or 2 vertices of $\{\infty_1,\infty_u\}$ in the construction above, the pairwise intersections can be of any of these magnitudes, but they will always intersect in the same vertices.

Theorem \ref{mt construction} provides sufficient conditions for the existence of a \UMSF$(v,m,k)$ if 
$m$ is a divisor of $v$. From that, we easily obtain the following corollary. 
\begin{corollary}\label{scattered_subsystems}
  Let $u\in\{1,2\}$, $k\geq 3$, and let $k\mid  m_j \mid  (2m+u)$ for each $1 \leq j \leq n$.
  Then there exists a \UMSF$(2ms+u, [m_j],k)$ whenever the following conditions hold:
  \begin{enumerate}
  \item $m_j$ is even or $(2m+u)/m_j$ is odd,
  \item $2n\leq 2s$ and $2s \equiv 2 \pmod{k}$,
  \end{enumerate}
  except when
  $(m_j, k)\in \{(6,3), (12,3)\}$, and except possibly when
  $(2m+u, m_j, k) = (12, 6, 6)$, or 
  $k$ is odd and $2m+u = 2m_j$,
  for some $j\in\{1,\ldots, n\}$.
\end{corollary}

Note that for values of the triple $(2m+u, m_j, k)$ determining a possible exception in Corollary \ref{scattered_subsystems}
it is possible for a \UMSF$(2ms+u, [m_j],k)$ to exist. However, our 
method cannot construct them 
because the \UMSF$(2m+u, m_j, k)$ to use in the construction
does not exist. It is possible that other construction methods would build a \UMSF$(2ms+u, [m_j],k)$.

\subsection{Nested subsystems}

A scenario complementary to the susbsystems being all disjoint is when the subsystems are completely nested, on vertex sets $M_1 \supset M_2 \supset \cdots \supset M_{n-1}$.  We modify our notation slightly for this section to make it less cumbersome in this specific context.

%  \begin{definition} \label{def_nested}
%  Let $v=m_0 > m_1 > \cdots > m_{n-1} > m_n=0$ be positive integers and $V=M_0 \supset  M_1 \supset \cdots \supset M_{n-1} \supset M_n=\varnothing$ be nested sets with $|M_j| = m_j$.  
%  For $1 \leq i \leq \lfloor \frac{m_0-1}{2} \rfloor$ and $0 \leq j < n$, 
%  let $U_{i,j}$ be a (possibly empty) $2$-regular graph such that
%  \begin{enumerate}
%    \item $V(U_{i,j}) = M_{j}\setminus M_{j+1}$ 
%    for $1 \leq i \leq \lfloor \frac{m_{j+1}-1}{2} \rfloor $;
%    \item $V(U_{i,j}) = M_{j}$  
%     for $\lfloor \frac{m_{j+1}+1}{2} \rfloor \leq i \leq 
%     \lfloor \frac{m_{j}-1}{2} \rfloor$;
%    \item $V(U_{i,j}) = \varnothing$     
%    for $\lfloor \frac{m_{j}+1}{2} \rfloor \leq i < 
%     \lfloor m_{0}/2 \rfloor - 1$.  
%  \end{enumerate}
%A {\em nested minisymposium factorization} 
%nMSF$(\{U_{i,j}:\;1 \leq i \leq 
%\lfloor \frac{m_0-1}{2} \rfloor,\; 0 \leq j < n\})$ is
%  \begin{itemize}
%  \item a factorization of $K_{v}^*$ with vertex set $V=M_0$ into 2-factors $F_i = \cup_{0 \leq j < n} U_{i,j}$;
%  \item for each $0 \leq \ell < n$, this factorization restricted to vertex set $M_{\ell}$ is a factorization of a graph isomorphic to $K_{m_{\ell}}^*$ into 2-factors 
%  $F_i' = \cup_{\ell \leq j < n} U_{i,j}$, 
%  for $1 \leq i < \lfloor \frac{m_{\ell}-1}{2} \rfloor$.
%    \end{itemize}
%  \end{definition}
  \begin{definition} \label{def_nested}
  Let $v=m_0 > m_1 > \cdots > m_{n-1} > m_n=0$ be non-negative integers.
  For $1 \leq i \leq \lfloor \frac{m_0-1}{2} \rfloor$ and $0 \leq j < n$, let $U_{i,j}$ be a $2$-regular graph of order
  \[
    |V(U_{i,j})| = 
    \begin{cases}
    m_j-m_{j+1}, & 
    \text{if $1 \leq i \leq \lfloor \frac{m_{j+1}-1}{2} \rfloor $,}\\
    m_j,  & 
    \text{if $\lfloor \frac{m_{j+1}+1}{2} \rfloor \leq i \leq 
     \lfloor \frac{m_{j}-1}{2} \rfloor$,}\\
    0,  & 
    \text{if $\lfloor \frac{m_{j}+1}{2} \rfloor \leq i \leq
     \lfloor \frac{m_0-1}{2} \rfloor$.} 
    \end{cases}
  \]
A {\em nested minisymposium factorization} 
nMSF$(\{U_{i,j}:\;1 \leq i \leq 
\lfloor \frac{m_0-1}{2} \rfloor,\; 0 \leq j < n\})$ is
a $2$-factorization $\mathcal{F}=\{F_i:  1 \leq i \leq
     \lfloor \frac{m_0-1}{2} \rfloor \}$ of $K_{v}^*$ such that
  \begin{itemize}
  \item $V(K_v^*)=M_0 \supset  M_1 \supset \cdots \supset M_{n-1} \supset M_n=\varnothing$ are nested sets with $|M_j| = m_j$;
  \item 
%for every $1 \leq i \leq
%     \lfloor \frac{m_0-1}{2} \rfloor$,
  each $2$-factor
  $F_i = \bigcup_{j=0}^{n-1} F_{i,j}$, where
  \[V(F_{i,j}) =
   \begin{cases}
     M_j\setminus M_{j+1} & \text{if $1 \leq i \leq \lfloor \frac{m_{j+1}-1}{2} \rfloor$}\\
     M_j & \text{if $\lfloor \frac{m_{j+1}+1}{2} \rfloor \leq i \leq \lfloor \frac{m_{j}-1}{2} \rfloor$}\\
     \varnothing & \text{otherwise},
   \end{cases}  
  \]
  and each $F_{i,j}$ is isomorphic to $U_{i,j}$;
  \item for every $0 \leq \ell < n$, 
  $\{\bigcup_{j=\ell}^{n-1} F_{i,j}:  1 \leq i \leq
     \lfloor \frac{m_\ell-1}{2} \rfloor \}$
is a $2$-factorization of a graph isomorphic to $K_{m_{\ell}}^*$.
    \end{itemize}
  \end{definition}
In other words, the factorization $\mathcal{F}$ of $K_v^*$
restricted to vertex set $M_{\ell}$ factorizes a graph isomorphic to $K_{m_{\ell}}^*$ into 2-factors whose structure is determined by the $U_{i,j}$s.  
  
  Our construction of nested minisymposium factorizations is most tidily expressed by defining holey factorizations.
\begin{definition} \label{def_hole}
  Given positive integers $v$ and $m$ with $v\geq m$, let 
  $$\cU = \left\{U_i:\; 1 \leq i \leq \left\lfloor\frac{v-1}{2}\right\rfloor \right\},$$ 
  be a collection of $2$-regular graphs on $v-m$ vertices for 
  $1 \leq i \leq \lfloor \frac{m-1}{2}\rfloor$, and
  on $v$ vertices for 
  $\lfloor \frac{m+1}{2}\rfloor \leq i \leq \lfloor \frac{v-1}{2}\rfloor$. 
  A {\em holey factorization} HF$(\cU)$ is a decomposition 
  $\mathcal{F}=\{F_i:\; 1 \leq i \leq \left\lfloor\frac{v-1}{2}\right\rfloor\}$ of
  $K_v^*-G$ (i.e., $K_v^*$ minus the edges of $G$) where
  each $F_i \cong U_i$ and $G \cong K_m^*$.
\end{definition}
If $v$ is even, then there is a 1-factor, $I_v$, on the vertices of $K_v^*$ whose edges are not present in $K_v^* - G$. If $v$ and $m$ are both even there is a 1-factor, $J_m$, on the vertices of $G$ which is a subgraph of $I_v$.  If $v$ is even and $m$ odd, then no edges of $I_v$ are induced on the vertices of $G$.  If $v$ is odd and $m$ is even, then there is a 1-factor $J_m$ on the vertices of $G$ whose edges are present in $K_v^*- G$.

By removing the 2-factors of a subsystem or ``filling the hole'' with them (making the $J_m$ in the hole coincide with the $I_m$ of the subsystem as required by the parities of $v$ and $m$) we have an equivalence between the existence of minisymposium factorizations and holey factorizations.
\begin{theorem}\label{thm_hole_eq}
  Let $\cT = \{T_i:\;1 \leq i \leq \lfloor \frac{m-1}{2} \rfloor\}$ be a $2$-factorization of $K_m^*$ and 
  $$\cU = \left\{U_i:\; 1 \leq i \leq \left\lfloor\frac{v-1}{2}\right\rfloor \right\},$$ 
be a collection of $2$-regular graphs on $v-m$ vertices for 
$1 \leq i \leq \lfloor \frac{m-1}{2} \rfloor$ and on $v$ vertices for $\lfloor \frac{m+1}{2}\rfloor \leq i \leq \lfloor \frac{v-1}{2}\rfloor $.  Then a HF$(\cU)$ exists if and only if a
\[\textstyle{
  \mbox{MSF}
  (\{U_i:\; \lfloor \frac{m+1}{2}\rfloor \leq i \leq \lfloor \frac{v-1}{2}\rfloor\},
   \{(T_i,U_i):\; 1 \leq i \leq \lfloor \frac{m-1}{2} \rfloor\})
   }
\]
exists.
\end{theorem}

Because in a nested minisymposium factorization the holes are nested and emptying or filling them does not affect the edges outside the hole, this equivalence extends to nested minisymposium factorizations and shows that they can be constructed exactly when the various holey factorizations of $K_{m_j}^*$ with holes of size $m_{j+1}$ exist.

\begin{theorem}
  \label{nested_holes}  
  Let $v=m_0 > m_1 > \cdots > m_{n-1} > m_n=0$ be positive integers.
  For $1 \leq i \leq \lfloor \frac{m_0-1}{2} \rfloor$ and $0 \leq j < n$, let $U_{i,j}$ be a $2$-regular graph of order
  \[
    |V(U_{i,j})| = 
    \begin{cases}
    m_j-m_{j+1}, & 
    \text{if $1 \leq i \leq \lfloor \frac{m_{j+1}-1}{2} \rfloor $,}\\
    m_j,  & 
    \text{if $\lfloor \frac{m_{j+1}+1}{2} \rfloor \leq i \leq 
     \lfloor \frac{m_{j}-1}{2} \rfloor$,}\\
    0,  & 
    \text{if $\lfloor \frac{m_{j}+1}{2} \rfloor \leq i \leq
     \lfloor \frac{m_0-1}{2} \rfloor$.} 
    \end{cases}
  \]
A nested minisymposium factorization
nMSF$(\{U_{i,j}:\;1 \leq i \leq \lfloor \frac{m_0-1}{2} \rfloor,\; 0 \leq j < n\})$ exists if and only if a HF$(\{U_{i,j}:\; 1 \leq i \leq \lfloor \frac{m_j-1}{2} \rfloor\})$ exists for each $0 \leq j < n$.
  \end{theorem}
  \begin{proof}
  The forward direction is proved simply by restricting the system to $M_j$ and removing the subsystem on $M_{j+1}$.  The converse is proved by a recursive construction starting with $j=n-1$: in this case,
  a HF$(\{U_{i,n-1}:\; 1 \leq i \leq \lfloor \frac{m_{n-1}-1}{2} \rfloor\})$ is an nMSF$(\{U_{i,n-1}:\;1 \leq i \leq \lfloor \frac{m_{n-1}-1}{2} \rfloor\})$, say $\mathcal{F}_{n-1}$. 
  
  At stage $j< n-1$,
use Theorem~\ref{thm_hole_eq} to construct 
 an nMSF$(\{U_{i,\ell}:\;1 \leq i \leq 
 \lfloor \frac{m_{j}-1}{2} \rfloor,\; j \leq \ell < n\})$, say $\mathcal{F}_j$,  by filling the hole in the 
 HF$(\{U_{i,j}:\; 1 \leq i \leq \lfloor \frac{m_{j}-1}{2} \rfloor\})$ with the 
 nMSF$(\{U_{i,\ell}:\;1 \leq i \leq 
 \lfloor \frac{m_{j+1}-1}{2} \rfloor,\; 
 j+1 \leq \ell < n\})$, denoted by $\mathcal{F}_{j+1}$, built at stage $j+1$.  
%
%  
%  This constructs an nMSF$(\{U_{i,\ell}:\;0 \leq i < \lfloor m_{j}/2 \rfloor,\; j \leq \ell < n\})$.
\end{proof}

Between the extremes of disjoint and nested subsystems,  there are factorizations with multiple subsystems with arbitrary intersections. Some structured instances of this much more general problem may be amenable to solution but we leave this to future work.

\section{Conclusions and Further Work}

We have introduced the minisymposium problem: a subsystem variant of the generalized Oberwolfach problem. This variant asks for a solution to a generalized Oberwolfach problem that contains a subsystem of a given size. 
When $v$, the number of vertices, is even, it is traditional in 2-factor decomposition problems 
to
ask for decompositions of $K_v-I$ where $I$ is a 1-factor. When the number of vertices in the system and the subsystem are both even, then we require that the 1-factor in the subsystem be a subgraph of the 1-factor in the full system. Therefore when the parities of the system and the subsystem agree, the problem becomes more tractable. When the parities are opposite, either the 1-factor of the full system must avoid the subsystem, or the edges of the 1-factor in the subsystem must be in 2-factors of the whole system.

Clearly, this is a very broad statement and we identify some particularly interesting cases, Hamiltonian* and uniform. In the Hamiltonian* minisymposium problem there are as few cycles as possible and in the uniform minisymposium problem all cycles are of the same length. 
We have shown that the work of Hilton and Johnson provides a complete solution for the Hamiltonian* minisymposium problem in Corollary~\ref{HMS exist}. In the case when $v \equiv m \equiv 2 \pmod 4$, Theorem~\ref{Hagg Result} uses this Hamiltonian* construction to provide a wide range of solutions when the resulting factors are all bipartite. In particular, a uniform factorization with $k\equiv 2\pmod 4$, $v\geq 2m$ and $v\equiv m\equiv k \pmod{2k}$ always exists. Corollary~\ref{recurse bipartite} can be used to extend this to uniform factorizations where $k \equiv 0 \pmod 4$ or $v \equiv m \equiv 0 \pmod {2k}$.

In Section~\ref{Uniform Factors} we considered the uniform case. We have solved a large part of the spectrum. In particular, 
when $k$ is
even or $v$ is odd, Corollary~\ref{m|v} gives all cases when $m\mid v$ and Corollary~\ref{k=m} completely solves all cases when $k=m$ 
has
the same parity as $v$.  Theorem~\ref{recursive} gives a powerful recursive construction which is applicable in cases where $m$ does not divide $v$. By applying it to the case when $k=3$, we obtain uniform factorizations with cycle lengths divisible by 3.  The case when $m$ is odd and $v$ is even seems to be the hardest.  Even in the simplest case when $k=3$, which has been well studied otherwise \cite{DRS_03_2, DRS_03_1, DRS_08, GR, Rees_Stinson_1, Rees_Stinson_2, Stinson, TRS}, the case with $v$ and $m$ having opposite parities has not been previously considered and remains open. 

While the Hamiltonian* problem is solved and we have made significant inroads into the uniform case, the general problem remains wide open. 
We expect that when $m$ and $v$ have the same parity solutions will be easier to find. When $m$ and $v$ have opposite parity we expect that $v$ odd with $m$ even is more tractable than the reverse.  Considering 2-factorizations where a solution to the Oberwolfach problem is known might be a good starting point.
A natural case to consider is the case when all factors are isomorphic i.e.\ $F_i\cong T_j \cup U_j$ for all $i$ and $j$.
Uniform factorizations are an example of this, but other variations are possible, for example, requiring all factors to be isomorphic to $C_{v-m} \cup C_m$.  Indeed, Theorem~\ref{Hagg Result} solves all these cases when the factors are bipartite and $v \equiv m \equiv 2 \pmod 4$, but this broader variant remains open.

More complex variants can also be considered.  We have briefly considered systems with multiple subsystems.  When these subsystems are completely nested the problem essentially reduces to the existence of the necessary ingredients as described in Theorem~\ref{nested_holes}.  Let $\{m_j\}$ be a multiset of subsystem sizes. When the subsystems are pairwise disjoint, $v$ is divisible by each $m_j$ and either $v$ is odd or at least one subsystem is even, then Lemma~\ref{disjoint_holes_minimal} and Corollary~\ref{disjoint_holes_thm} use Theorem~\ref{Liu} to construct a \UMSF$(v,\{m_j\},k)$ for all but a finite number of admissible $v$. Even in the seemingly simple case when the subsystems are all disjoint the problem remains generally open even for the uniform case.  Further partial results are obtained when the subsystems are scattered, that is, when no two minisymposia have meetings taking place on the same day. Cycle frames in Theorem \ref{CF} allow us to construct uniform factorizations as described in Lemma~\ref{scattered_subsystems_minimal} and Theorem~\ref{scattered_subsystems}.  The more general case when the intersections of multiple subsystems are arbitrary is completely open.

\section*{Acknowledgements}
We thank the anonymous referees for their many useful suggestions that helped strongly improve this paper. 

P.\ Danziger has received support from NSERC Discovery Grants 
RGPIN-2022-03816. B.\ Stevens recieved support from NSERC Discovery Grant RGPIN-2017-06392. 
T.\ Traetta has received support from GNSAGA of Istituto Nazionale di Alta Matematica.

\bibliographystyle{abbrv}
\bibliography{refs.bib}

\begin{thebibliography}{10}

\bibitem{Alspach_Haggkvist_85}
B.~Alspach and R.~H\"{a}ggkvist.
\newblock Some observations on the {O}berwolfach problem.
\newblock {\em J. Graph Theory}, 9(1):177--187, 1985.

\bibitem{ASSW}
B.~Alspach, P.~J. Schellenberg, D.~R. Stinson, and D.~Wagner.
\newblock The {O}berwolfach problem and factors of uniform odd length cycles.
\newblock {\em J. Combin. Theory Ser. A}, 52(1):20--43, 1989.

\bibitem{AsplundEtAl}
J.~Asplund, D.~Kamin, M.~Keranen, A.~Pastine, and S.~\"{O}zkan.
\newblock On the {H}amilton-{W}aterloo problem with triangle factors and
  {$C_{3x}$}-factors.
\newblock {\em Australas. J. Combin.}, 64:458--474, 2016.

\bibitem{BS81}
Z.~Baranyai and G.~R. Sz\'{a}sz.
\newblock Hamiltonian decomposition of lexicographic product.
\newblock {\em J. Combin. Theory Ser. B}, 31(3):253--261, 1981.

\bibitem{BonviciniBuratti}
S.~Bonvicini and M.~Buratti.
\newblock Octahedral, dicyclic and special linear solutions of some
  {H}amilton-{W}aterloo problems.
\newblock {\em Ars Math. Contemp.}, 14(1):1--14, 2018.

\bibitem{BryantDanziger}
D.~Bryant and P.~Danziger.
\newblock On bipartite 2-factorizations of {$K_n-I$} and the {O}berwolfach
  problem.
\newblock {\em J. Graph Theory}, 68(1):22--37, 2011.

\bibitem{BryantDanzigerDean}
D.~Bryant, P.~Danziger, and M.~Dean.
\newblock On the {H}amilton-{W}aterloo problem for bipartite 2-factors.
\newblock {\em J. Combin. Des.}, 21(2):60--80, 2013.

\bibitem{BCDT17}
M.~Buratti, H.~Cao, D.~Dai, and T.~Traetta.
\newblock A complete solution to the existence of {$(k,\lambda)$}-cycle frames
  of type {$g^u$}.
\newblock {\em J. Combin. Des.}, 25(5):197--230, 2017.

\bibitem{BZ}
M.~Buratti and F.~Zuanni.
\newblock {P}erfect {C}ayley {D}esigns as {G}eneralizations of {P}erfect
  {M}endelsohn {D}esigns.
\newblock {\em Des. Codes Cryptogr.}, 23:233–--248, 2001.

\bibitem{BDT1}
A.~C. Burgess, P.~Danziger, and T.~Traetta.
\newblock On the {H}amilton-{W}aterloo problem with odd orders.
\newblock {\em J. Combin. Des.}, 25(6):258--287, 2017.

\bibitem{BDT3}
A.~C. Burgess, P.~Danziger, and T.~Traetta.
\newblock On the {H}amilton-{W}aterloo problem with cycle lengths of distinct
  parities.
\newblock {\em Discrete Math.}, 341(6):1636--1644, 2018.

\bibitem{BDT2}
A.~C. Burgess, P.~Danziger, and T.~Traetta.
\newblock On the {H}amilton-{W}aterloo problem with odd cycle lengths.
\newblock {\em J. Combin. Des.}, 26(2):51--83, 2018.

\bibitem{CaElKhoVan04}
N.~J. Cavenagh, S.~I. El-Zanati, A.~Khodkar, and C.~Vanden~Eynden.
\newblock On a generalization of the {O}berwolfach problem.
\newblock {\em J. Combin. Theory Ser. A}, 106(2):255--275, 2004.

\bibitem{Handbook}
C.~J. Colbourn and J.~H. Dinitz, editors.
\newblock {\em Handbook of combinatorial designs}.
\newblock Discrete Mathematics and its Applications (Boca Raton). Chapman \&
  Hall/CRC, Boca Raton, FL, second edition, 2007.

\bibitem{Costa20}
S.~Costa.
\newblock A complete solution to the infinite {O}berwolfach problem.
\newblock {\em J. Combin. Des.}, 28(5):366--383, 2020.

\bibitem{hamiltonian_paper}
P.~Danziger, E.~Mendelsohn, B.~Stevens, and T.~Traetta.
\newblock Hamiltonian* minisymposium factorizations.
\newblock preprint (2023).

\bibitem{DanzigerQuattrocchiStevens}
P.~Danziger, G.~Quattrocchi, and B.~Stevens.
\newblock The {H}amilton-{W}aterloo problem for cycle sizes 3 and 4.
\newblock {\em J. Combin. Des.}, 17(4):342--352, 2009.

\bibitem{DRS_03_2}
D.~Deng, R.~Rees, and H.~Shen.
\newblock Further results on nearly {K}irkman triple systems with subsystems.
\newblock {\em Discrete Math.}, 270(1-3):99--114, 2003.

\bibitem{DRS_03_1}
D.~Deng, R.~Rees, and H.~Shen.
\newblock On the existence and application of incomplete nearly {K}irkman
  triple systems with a hole of size 6 or 12.
\newblock {\em Discrete Math.}, 261(1-3):209--233, 2003.
\newblock Papers on the occasion of the 65th birthday of Alex Rosa.

\bibitem{DRS_08}
D.~Deng, R.~Rees, and H.~Shen.
\newblock On the existence of nearly {K}irkman triple systems with subsystems.
\newblock {\em Des. Codes Cryptogr.}, 48(1):17--33, 2008.

\bibitem{ElTipVan02}
S.~I. El-Zanati, S.~K. Tipnis, and C.~Vanden~Eynden.
\newblock A generalization of the {O}berwolfach problem.
\newblock {\em J. Graph Theory}, 41(2):151--161, 2002.

\bibitem{GR}
G.~Ge and R.~Rees.
\newblock On group-divisible designs with block size four and group-type
  {$6^um^1$}.
\newblock {\em Discrete Math.}, 279(1-3):247--265, 2004.
\newblock In honour of Zhu Lie.

\bibitem{MR4012874}
S.~Glock, F.~Joos, J.~Kim, D.~K\"{u}hn, and D.~Osthus.
\newblock Resolution of the {O}berwolfach problem.
\newblock {\em Acta Math. Univ. Comenian. (N.S.)}, 88(3):735--741, 2019.

\bibitem{Hag}
R.~H\"{a}ggkvist.
\newblock A lemma on cycle decompositions.
\newblock In {\em Cycles in graphs ({B}urnaby, {B}.{C}., 1982)}, volume 115 of
  {\em North-Holland Math. Stud.}, pages 227--232. North-Holland, Amsterdam,
  1985.

\bibitem{MR1865547}
A.~J.~W. Hilton and M.~Johnson.
\newblock Some results on the {O}berwolfach problem.
\newblock {\em J. London Math. Soc. (2)}, 64(3):513--522, 2001.

\bibitem{Hoffman_Schellenberg_91}
D.~G. Hoffman and P.~J. Schellenberg.
\newblock The existence of {$C_k$}-factorizations of {$K_{2n}-F$}.
\newblock {\em Discrete Math.}, 97(1-3):243--250, 1991.

\bibitem{KeranenOzkan}
M.~S. Keranen and S.~\"{O}zkan.
\newblock The {H}amilton-{W}aterloo problem with 4-cycles and a single factor
  of {$n$}-cycles.
\newblock {\em Graphs Combin.}, 29(6):1827--1837, 2013.

\bibitem{KeranenPastine}
M.~S. Keranen and A.~Pastine.
\newblock A generalization of the {H}amilton-{W}aterloo problem on complete
  equipartite graphs.
\newblock {\em J. Combin. Des.}, 25(10):431--468, 2017.

\bibitem{LeiShen}
H.~Lei and H.~Shen.
\newblock The {H}amilton-{W}aterloo problem for {H}amilton cycles and
  triangle-factors.
\newblock {\em J. Combin. Des.}, 20(7):305--316, 2012.

\bibitem{Liu00}
J.~Liu.
\newblock A generalization of the {O}berwolfach problem and
  {$C_t$}-factorizations of complete equipartite graphs.
\newblock {\em J. Combin. Des.}, 8(1):42--49, 2000.

\bibitem{Liu03}
J.~Liu.
\newblock The equipartite {O}berwolfach problem with uniform tables.
\newblock {\em J. Combin. Theory Ser. A}, 101(1):20--34, 2003.

\bibitem{OdabasiOzkan}
U.~Odaba\c{s}\i and S.~\"{O}zkan.
\newblock The {H}amilton-{W}aterloo problem with {$C_4$} and {$C_m$} factors.
\newblock {\em Discrete Math.}, 339(1):263--269, 2016.

\bibitem{Piotrowski}
W.-L. Piotrowski.
\newblock The solution of the bipartite analogue of the {O}berwolfach problem.
\newblock {\em Discrete Math.}, 97(1-3):339--356, 1991.

\bibitem{Rees_Stinson_1}
R.~Rees and D.~R. Stinson.
\newblock Kirkman triple systems with maximum subsystems.
\newblock {\em Ars Combin.}, 25:125--132, 1988.

\bibitem{Rees_Stinson_2}
R.~Rees and D.~R. Stinson.
\newblock On the existence of {K}irkman triple systems containing {K}irkman
  subsystems.
\newblock {\em Ars Combin.}, 26:3--16, 1988.

\bibitem{Stinson}
D.~R. Stinson.
\newblock Frames for {K}irkman triple systems.
\newblock {\em Discrete Math.}, 65(3):289--300, 1987.

\bibitem{TRS}
S.~Tang and H.~Shen.
\newblock Embeddings of nearly {K}irkman triple systems.
\newblock {\em J. Statist. Plann. Inference}, 94(2):327--333, 2001.
\newblock Second Shanghai Conference on Designs, Codes and Finite Geometries
  (1996).

\bibitem{Traetta_13}
T.~Traetta.
\newblock A complete solution to the two-table {O}berwolfach problems.
\newblock {\em J. Combin. Theory Ser. A}, 120(5):984--997, 2013.

\bibitem{WangCao}
L.~Wang and H.~Cao.
\newblock A note on the {H}amilton-{W}aterloo problem with {$C_8$}-factors and
  {$C_m$}-factors.
\newblock {\em Discrete Math.}, 341(1):67--73, 2018.

\bibitem{WangChenCao}
L.~Wang, F.~Chen, and H.~Cao.
\newblock The {H}amilton-{W}aterloo problem for {$C_3$}-factors and
  {$C_n$}-factors.
\newblock {\em J. Combin. Des.}, 25(9):385--418, 2017.

\bibitem{West}
D.~B. West.
\newblock {\em Introduction to graph theory}.
\newblock Prentice Hall, Inc., Upper Saddle River, NJ, 2001.

\end{thebibliography}
\end{document}